\documentclass[12pt]{amsart}

\usepackage{amsmath,amssymb,amsthm}
 \usepackage{setspace}
\theoremstyle{plain}

\theoremstyle{plain}
\theoremstyle{plain}
\newtheorem{Theo}{Theorem}[section]

\newtheorem{lem}[Theo]{Lemma}

 \newtheorem{prop}[Theo]{Proposition}

\newcommand{\ZZ}{\mathbb{Z}}  
  
\newcommand{\NN}{\mathbb{N}}

\newcommand{\RR}{\mathbb{R}}

\usepackage{amsfonts}
 
\usepackage{fancyhdr}
\usepackage{amsfonts,graphicx}

 \newtheorem{Rema}{Remark }

\newcommand{\QG}{(\textnormal{QG}_{\alpha})}

\begin{document}

\title[On the global solutions of the super-critical $2$D Q-G equation in  Besov spaces]{On the global solutions of the super-critical $2$D quasi-geostrophic equation in  Besov spaces}

\author[T. Hmidi]{Taoufik Hmidi}
\address{IRMAR, Universit\'e de Rennes 1\\ Campus de
Beaulieu\\ 35~042 Rennes cedex\\ France}
\email{thmidi@univ-renne1.fr}
\author[S. Keraani]{Sahbi Keraani}
\address{IRMAR, Universit\'e de Rennes 1\\ Campus de
Beaulieu\\ 35~042 Rennes cedex\\ France}
\email{sahbi.keraani@univ-rennes1.fr}


\begin{abstract}
In this paper we study the super-critical $2$D dissipative quasi-geostrophic equation.
 We obtain some regularization effects allowing us to prove   global well-posedness result 
  for small initial data  lying in critical Besov spaces constructed over Lebesgue spaces $L^p,$ \mbox{with $p\in[1,\infty].$ Local results for arbitrary initial data are also given.}
\end{abstract}

\maketitle

\maketitle
\section{Introduction}
This paper deals with the Cauchy problem for the two-dimensional dissipative quasi-geostrophic equation
$$
\QG\left\lbrace
\begin{array}{l}\partial_{t}\theta+v\cdot\nabla\theta+|\hbox{D}|^{{\alpha}} \theta=0\\
\theta_{\vert t =0}=\theta^0,
\end{array}
\right.$$
where  the scalar function $\theta$ represents the potential temperature and $\alpha \in[0,2].$ The velocity $v=(v^1,v^2)$ is determined by $\theta$ through a stream function $\psi$, namely
  $$
  v=(-\partial_{2}\psi,\partial_{1}\psi), \quad\textnormal{with}\quad|\hbox{D}|\psi=\theta.
  $$
Here, the differential  operator $|\hbox{D}|=\sqrt{-\Delta}$ is defined in a standard fashion  through its Fourier transform: $\mathcal{F}(|D|u)=|\xi|\mathcal{F}u.$  The above relations can be rewritten  as
$$
v=(-\partial_{2}|\hbox{D}|^{-1}\theta,\partial_{1}|\hbox{D}|^{-1}\theta)=(-{R}_{2}\theta,{R}_{1}\theta),
$$
 where ${R}_{i} (i=1,2)$ are Riesz transforms.
 
 First we notice that solutions for $\QG$ equation are scaling invariant in the following sense: if $\theta$ is a solution and $\lambda>0$ then $\theta_{\lambda}(t,x)=\lambda^{\alpha-1}\theta(\lambda^\alpha t,\lambda x)$ is also a solution of $\QG$ equation. From the definition of the homogeneous Besov spaces, described in next section, one  can show that the norm of $\theta_{\lambda}$ in the space  $\dot B_{p,r}^{1+\frac{2}{p}-\alpha},$ with $p,\,r\in[1,\infty],$ is quasi-invariant. That is, there exists a pure constant $C>0$ such that for\mbox { every $\lambda, t>0$}
 $$
 C^{-1}\|\theta_{\lambda}(t)\|_{\dot B_{p,r}^{1+\frac{2}{p}-\alpha}}\leq\|\theta(\lambda^\alpha t)\|_{\dot B_{p,r}^{1+\frac{2}{p}-\alpha}}\leq C\|\theta_{\lambda}(t)\|_{\dot B_{p,r}^{1+\frac{2}{p}-\alpha}}.
  $$
 
  Besides its intrinsic mathematical importance the $\QG$
     equation serves as  a $2$D models arising 
  in geophysical fluid dynamics, for more details about the subject \mbox{see \cite{C-M-T,Ped}} and the references therein.   
Recently the $\QG$ equation has been intensively investigated and  much attention is carried to the problem of global existence.  For the sub-critical case
   $(\alpha>1)$ the theory seems to be in a satisfactory state. Indeed, global existence and uniqueness for arbitrary  initial data  are established in various function spaces \mbox{(see for example  \cite{C-W,Res}).} However in the critical case, that is $\alpha=1,$ Constantin \mbox{et {\it al.} \cite{C-C-W}} showed  the global existence in Sobolev space $H^1$   under smallness  assumption of the $L^\infty$-norm of the initial temperature $\theta^0$ but the uniqueness is proved for initial data in $H^2.$ Many other relevant results can be found in \cite{C-C,Ju,Ju1}. The super-critical case $\alpha<1$ seems  harder to deal with and work on this subject has just started to appear. In \cite{Cha} the global existence and uniqueness are established for data in critical Besov space $B_{2,1}^{2-\alpha}$ with a small $\dot{B}_{2,1}^{2-\alpha}$ norm. This result was  improved by N. Ju \cite{Ning1}  for small initial data in $H^s$ with $s\geq2-\alpha.$ We would like to point out that all these spaces are constructed over Lebesgue space $L^2$ and the same problem for general Besov space $B_{p,r}^s$ is not yet well explored and few results are obtained in this subject. In \cite{Wu1}, Wu  proved the  global existence and uniqueness for small initial data in $C^r\cap L^q$ with $r>1$ and $q\in]1,\infty[,$ which is not a scaling space. We can also mention the paper \cite{Wu2} in which global well-posedness is established for small initial data in $B_{2,\infty}^s\cap B_{p,\infty}^s,$ with $s>2-\alpha$ and  $p=2^N.$ 
   
   The main goal of the present paper is to study existence and uniqueness problems in the super-critical case   when initial data belong to  inhomogeneous critical Besov \mbox{spaces  $B_{p,1}^{1+\frac{2}{p}-\alpha},$} with $p\in[1,\infty].$  
   
   Our first main result reads as follows. 
\begin{Theo}\label{Thm1} Let $\alpha\in[0,1[,$ $p\in[1,\infty]$ and $s\geq s_{c}^p,$ with $ s_{c}^p=1+\frac{2}{p}-\alpha$ and define
 $$
 {\mathcal{X}}_{p}^s=\left\lbrace
\begin{array}{l}B_{p,1}^s,\quad\hbox{if}\quad p<\infty,\\
B_{\infty,1}^s\cap \dot{B}_{\infty,1}^0,\quad\hbox{otherwise}. \
\end{array}
\right. $$
Then for $\theta^0\in {\mathcal{X}}_{p}^s$ there exists $T>0$ such that the $\QG$ equation has a unique solution   $\theta $ belonging to ${\it C}([0,T];{\mathcal{X}}_{p}^s)\cap L^1([0,T];\dot{B}_{p,1}^{s+\alpha}).$

In addition, there exists an absolute constant $\eta>0$ such that if 
$$
\|\theta^0\|_{\dot{B}_{\infty,1}^{1-\alpha}}\leq\eta,
$$
 then one can take $T=+\infty$.
 
\end{Theo}
\begin{Rema} We observe that in our global existence result we make only a smallness assumption of the data in Besov space
 $\dot B_{\infty,1}^{1-\alpha}$ which contains the increasing Besov chain  spaces $\{\dot{B}_{p,1}^{s_{c}^p}\}_{p\in[1,\infty]}.$
\end{Rema}
\begin{Rema}
In the case of $s>s_{c}^p$ we have the following lower bound for the local time existence. There exists a nonnegative  constant $C$ such that
$$
T\geq {C}{\|\theta^0\|_{\dot{B}_{\infty,1}^{s-\frac{2}{p}}}^{\frac{-\alpha}{s-s_{c}^p}}}.
$$
However in the critical case $s=s^p_{c}$ the local time existence is bounded from below by
$$
\sup\Big\{t\geq 0,\,\sum_{q\in\ZZ}(1-e^{-ct2^{q\alpha}})^{\frac{1}{2}}2^{q(1-\alpha)}\|\Delta_{q}\theta^0\|_{L^\infty}\leq \eta \Big\},
$$
where $\eta$ is an absolute nonnegative constant.
\end{Rema}

The proof relies essentially on some new estimates for transport-diffusion equation

 $$
 (\textnormal{TD}_{\alpha})\left\lbrace
\begin{array}{l}
\partial_t \theta+v\cdot\nabla \theta+|\textnormal{D}|^{\alpha} \theta=f\\
{\theta}_{| t=0}=\theta^{0},\\
\end{array}
\right.
$$
where the unknown is the scalar function $\theta.$ Our second main result reads as follows
 \begin{Theo}\label{Thm3}

 Let $s\in]-1,1[,\,\alpha\in[0,1[,\,(p,r)\in[1,+\infty]^2,$ $f\in L^1_{\textnormal{loc}}(\RR_{+}; \dot{B}_{p,1}^s)$ and $v$ be a divergence free  vector field belonging to $L^1_{\textnormal{loc}}(\RR_{+};\textnormal{Lip}(\RR^d)).$ 
  We consider a smooth solution $\theta$ of the transport-diffusion equation $(\textnormal{TD}_{\alpha})$
, then there exists a constant $C$ depending only on $s$ and $\alpha$ such that
 $$
 \|\theta\|_{\widetilde{L^r_{t}}\dot{B}_{p,1}^{s+\frac{\alpha}{r}}}\leq C e^{C\int_{0}^t\|\nabla v(\tau)\|_{L^\infty}d\tau}\Big(  \|\theta^0\|_{\dot{B}_{p,1}^s}+\|f\|_{L^1_{t}\dot{B}_{p,1}^s}\Big).
 $$
 Besides if $v=\nabla^\bot|\textnormal{D}|^{-1}\theta$ then the above estimate is valid for all $s>-1.$
    \end{Theo}
We use for the proof a new approach  based on
 Lagrangian coordinates  combined with paradifferential calculus and a new commutator estimate. This idea has been recently used by the first author to treat the  two-dimenional Navier-Stokes vortex patches \cite{Hmidi}.
\begin{Rema}
The estimates of Theorem $\ref{Thm3}$ hold true for Besov spaces $\dot B_{p,m}^s,$ with $m\in[1,\infty].$ The proof can be done  strictly in  the same line as the case $m=1.$
It should be also mentioned that we can derive similar results for 
  inhomogeneous Besov spaces. 
  \end{Rema}
{\bf Notation:} Throughout the paper, $C$ stands for a  constant which may be different in each occurrence. We shall sometimes use the notation $A\lesssim B$ instead of $A\leq CB$ and $A\approx B$ means that $A\lesssim B$ and $B\lesssim A.$

The rest of this paper is structured as follows. In next section we recall some basic results on Littlewood-Paley theory and we give some useful lemmas. Section $3$ is devoted to the proof of a new commutator estimate while  sections $4$ and $5$ are dealing successively  with  the  proofs of Theorem \ref{Thm3} and \ref{Thm1}. We give in the end of this paper an appendix.
\section{Preliminaries}
In this  preparatory section, we provide the definition of some function spaces based on the so-called Littlewood-Paley decomposition and we  review some important lemmas that will 
be used  constantly in the following pages.

We start with the dyadic decomposition. Let $\varphi\in C^\infty_{0}(\RR^d)$ be supported in the ring $\mathcal{C}:=\{ \xi\in\RR^d,\frac{3}{4}\leq|\xi|\leq\frac{8}{3}\}$ and such that 
$$
\sum_{q\in\ZZ}\varphi(2^{-q}\xi)=1 \quad\hbox{for}\quad \xi\neq 0.
$$
We define also the function $\;\chi(\xi)=1-\sum_{q\in\NN}\varphi(2^{-q}\xi).$ Now  for  $u\in{\mathcal S}'$ we set
  $$
\Delta_{-1}u=\chi(\hbox{D})u;\, \forall
 q\in\NN,\;\Delta_qu=\varphi(2^{-q}\hbox{D})u\quad\hbox{ and  }\;\forall\,\,q\in\ZZ,\,\,\dot{\Delta}_{q}u=\varphi(2^{-q}\textnormal{D})u.$$ 
The following low-frequency cut-off will be also used:
$$
S_q u=\sum_{-1\leq j\leq q-1}\Delta_{j}u\quad\hbox{and}\quad \dot{S}_{q}u=\sum_{j\leq q-1}\dot\Delta_{j}u.
$$
We caution that we shall sometimes use the notation $\Delta_{q}$ instead of $\dot\Delta_{q}$ and this will be tacitly understood from the context. \\
Let us now recall the definition of  Besov spaces through dyadic decomposition. 

Let $(p,m)\in[1,+\infty]^2$ and $s\in\RR,$ then the inhomogeneous \mbox{space  $B_{p,m}^s$} is 
the set of tempered distribution $u$ such that
$$
\|u\|_{B_{p,m}^s}:=\Big( 2^{qs}
\|\Delta_q u\|_{L^{p}}\Big)_{\ell ^{m}}<\infty.
$$
To define the homogeneous Besov spaces we first denote by $\mathcal{S}'/\mathcal{P}$ the space of tempered distributions modulo polynomials.
Thus we define the space $\dot B_{p,r}^s$ as the set of distribution $u\in\mathcal{S}'/\mathcal{P}$ such that 
$$
\|u\|_{\dot{B}_{p,m}^s}:=\Big( 2^{qs}
\|\dot{\Delta}_q u\|_{L^{p}}\Big)_{\ell ^{m}}<\infty.
$$
We point out that if $s>0$ then we have $B_{p,m}^s=\dot B_{p,m}^s\cap L^{p}$ and
$$
\|u\|_{B_{p,m}^s}\approx\|u\|_{\dot B_{p,m}^s}+\|u\|_{L^{p}}.$$
Another characterization of homogeneous Besov spaces that will be needed later is given as follows (see \cite{Triebel}). For $s\in]0,1[, p,m\in[1,\infty]$
\begin{equation}
\label{equivalence}
C^{-1}\|u\|_{\dot{B}_{p,m}^s}\leq 
\Big(\int_{\RR^d}\frac{\|u(\cdot-x)-u(\cdot)\|^{m}_{L^p}}{|x|^{sm}}\frac{dx}{|x|^d}\Big)^{\frac{1}{m}}\leq C\|u\|_{\dot{B}_{p,m}^s},
\end{equation}
with the usual modification if $m=\infty$.\\
In our next study we require two kinds of coupled space-time Besov spaces. The first one is defined in the following manner: for  $T>0$ \mbox{and $m\geq1,$} we denote by $L^r_{T}\dot B_{p,m}^s$ the set of all tempered distribution $u$ satisfying
$$
\|u\|_{L^r_{T}\dot B_{p,r}^s}:= \Big\|\Big( 2^{qs}
\|\dot\Delta_q u\|_{L^{p}}\Big)_{\ell ^{m}}\Big\|_{L^r_{T}}<\infty.$$
The second mixed space is  $\widetilde L^r_{T}{\dot B_{p,m}^s}$ which is the set of  tempered distribution $u$ satisfying
 $$
 \|u\|_{ \widetilde L^r_{T}{\dot B_{p,m}^s}}:= \Big( 2^{qs}
\|\dot\Delta_q u\|_{L^r_{T}L^{p}}\Big)_{\ell ^{m}}<\infty .$$
We can define by the same way the spaces $L^r_{T} B_{p,m}^s$ and $\widetilde L^r_{T}{ B_{p,m}^s}.$ 

The following embeddings  are a direct consequence of Minkowski's inequality. 

Let $s\in\RR,$ $r\geq1$ and $\big(p,m\big)\in[1,\infty]^2,$ then we have 
\begin{eqnarray}\label{lemm4}
L^r_{T}\dot B_{p,m}^s&\hookrightarrow&\widetilde L^r_{T}{\dot B_{p,m}^s},\,\textnormal{if}\quad  m\geq r\quad\hbox{and}\\
\nonumber\widetilde L^r_{T}{\dot B_{p,m}^s}&\hookrightarrow& L^r_{T}\dot B_{p,m}^s,\, \textnormal{if}\quad 
r\geq m.
\end{eqnarray}
Another classical result that will be frequently used here is the so-called Bernstein inequalities (see \cite{Ch1} and the references therein): there exists $C$ such that for every function $u$ and for every $q\in\ZZ,$ we have
\begin{eqnarray*}
\sup_{|\alpha|=k}\|\partial ^{\alpha}S_{q}u\|_{L^b}&\leq& C^k\, 2^{q(k+d(\frac{1}{a}-\frac{1}{b}))}\|S_{q}u\|_{L^a},\quad\hbox{for}\quad b\geq a,\\
 C^{-k}2^{qk}\|\dot\Delta_{q}u\|_{L^a}&\leq&\sup_{|\alpha|=k}\|\partial ^{\alpha}\dot\Delta_{q}u\|_{L^a}\leq C^k2^{qk}\|\dot\Delta_{q}u\|_{L^a}.
\end{eqnarray*}
 It is worth pointing out that the above inequalities hold true if we replace the derivative $\partial^\alpha$ by fractional derivative $|\textnormal{D}|^\alpha.$
 According to Bernstein inequalities one can show the following embeddings
 $$
 \dot{B}_{p,m}^s\hookrightarrow\dot{B}_{p_{1},m_{1}}^{s-d(\frac{1}{p}-\frac{1}{p_{1}})},\quad\hbox{for}\quad p\leq p_{1}\quad\hbox{and}\quad m\leq m_{1}.
 $$
 
Now let us  we recall the following commutator lemma (see \cite{Ch1,rd23} and the references therein).
\begin{lem}\label{lemm12}
Let  $p,r\in[1,\infty], 1=\frac{1}{r}+\frac{1}{\bar r},  \rho_{1}<1, \rho_{2}<1$ and $v$ be a divergence free  vector field of $\RR^d.$ Assume in addition that
$$
\rho_{1}+\rho_{2}+d\min\{1,{2}/{p}\}>0\quad\hbox{and}\quad\rho_{1}+{d}/{p}>0.$$
Then we have
$$
\sum_{q\in\ZZ}2^{q(\frac{d}{p}+\rho_{1}+\rho_{2}-1)}\big\|[\dot \Delta_q, v\cdot \nabla ]u\big\|_{L^1_{t}L^{p}}\lesssim  \|
v\|_{\widetilde L^r_{t}\dot B_{p,1}^{\frac{d}{p}+\rho_{1}}}\|u \|_{\widetilde L^{\bar r}_{t}\dot B_{p,1}^{\frac{d}{p}+\rho_{2}}}.
$$ 
Moreover we have for $s\in]-1,1[$
$$
\sum_{q\in\ZZ}2^{qs}\big\|[\dot \Delta_q, v\cdot \nabla ]u\big\|_{L^{p}}\lesssim \|\nabla
v\|_{L^{\infty}}\|u\|_{\dot B_{p,1}^s}.$$
In addition this estimate holds true for all $s>-1$ if $v=\nabla^\perp|\textnormal{D}|^{-1}u.$
\end{lem}

 The following result describes the action of the semi-group operator $e^{t|\textnormal{D}|^\alpha}$ 
on distributions whose Fourier transform is supported in a ring.
\begin{prop}\label{l:5}
Let $\mathcal{C}$ be a ring and $\alpha\in\RR_{+}.$ There exists a positive constant $C$ such
 that for any $p\in[1;+\infty],$ for any couple $(t,\lambda)$ of positive real numbers, we have
$$
\textnormal{supp} \mathcal{ F}u\subset\lambda\mathcal{C}\Rightarrow 
\|e^{t|\textnormal{D}|^\alpha}u\|_{L^p}\leq Ce^{-C^{-1}t\lambda^{\alpha}}\|u\|_{L^p}.
$$
\end{prop}
\begin{proof}
We will imitate the same idea of \cite{Ch2}. Let $\phi\in \mathcal{D}(\RR^d\backslash\{ 0\}),$ radially and whose value is identically $1$ near the ring $\mathcal{C}.$ Then we have
\begin{eqnarray*}
e^{t|\textnormal{D}|^\alpha }u=\phi(\lambda^{-1}|\textnormal{D}|) u=h_{\lambda}*u,
\end{eqnarray*}
where
$$
h_{\lambda}(t,x)=\frac{1}{(2\pi)^d}\int_{\RR^d}\phi(\lambda^{-1}\xi)e^{-t\vert\xi\vert^\alpha}e^{i<x,\xi>}d\xi.
$$
We set 
$$ \bar h_{\lambda}(t,x):=\lambda^{-d}h(t,\lambda^{-1}x)=\frac{1}{(2\pi)^d}
\int_{\RR^d}\phi(\xi)e^{-t\lambda^\alpha\vert\xi\vert^\alpha}e^{i<x,\xi>}d\xi.
$$
Now to prove the proposition it suffices to show that $\|\bar h_{\lambda}(t)\|_{L^1}\leq Ce^{-C^{-1}t\lambda^{\alpha}}.$ 
For this purpose we  write with the aid of an integration by parts
$$
(1+\vert x\vert^2)^d{\bar h_{\lambda}}(x)=\frac{1}{(2\pi)^d}
\int_{\RR^d}(\textnormal{Id}-\Delta_{\xi})^d\big(\phi(\xi)e^{-t\lambda^\alpha\vert\xi\vert^\alpha}\big)e^{i<x,\xi>}d\xi.$$
From Leibnitz's formula, we have
$$
(\textnormal{Id}-\Delta_{\xi})^d\big(\phi(\xi)e^{-t\lambda^\alpha \xi^\alpha}\big)=\sum_{\vert
\gamma\vert\leq 2d\\
\atop \beta\leq\gamma}C_{\gamma,\beta}\partial^{\gamma-\beta}\phi(\xi)
\partial^\beta e^{-t\lambda^\alpha\xi^\alpha}.$$
As $\phi$ is supported in a ring that does not contain some neighbourhood of zero then we get for $\xi\in \hbox{supp }\phi$
\begin{eqnarray*}
\arrowvert\partial^\beta e^{-t\lambda^\alpha\vert\xi\vert^\alpha}
\arrowvert&\leq& C_{\beta}(1+t\lambda^\alpha)^{\vert\beta\vert}
e^{-t\lambda^\alpha\vert\xi\vert^\alpha}, \;\forall\, \xi\in \hbox{ supp }\phi\\
&\leq& C_{\beta}e^{-C^{-1}t\lambda^\alpha}.
\end{eqnarray*}
Thus we find that
$$
\big{|}(\textnormal{Id}-\Delta_{\xi})^d\big(\phi(\xi)e^{-t\lambda^\alpha \xi^\alpha}\big)\big{|}\leq  Ce^{-C^{-1}t\lambda^\alpha}\sum_{\vert
\gamma\vert\leq 2d\\
\atop \beta\leq\gamma}C_{\gamma,\beta}|\partial^{\gamma-\beta}\phi(\xi)|.
$$
Since the term of the right-hand side belongs to $L^1(\RR^d),$ then we deduce that
$$
(1+\vert x\vert^2)^d\arrowvert{\bar h_{\lambda}}(x)\arrowvert\leq Ce^{-C^{-1}t\lambda^\alpha}.$$
This completes the proof of the proposition.
\end{proof}

\section{Commutator estimate}
The main result of this section is the following estimate that will play a  crucial role for the proof of Theorem \ref{Thm3}.
\begin{prop}\label{pr:1}
Let
$f\in \dot{B}_{p,1}^\alpha$ with $\alpha\in [0,1[$ and $p\in[1,+\infty],$ and let $\psi$ be a  Lipshitz  measure-preserving homeomorphism \mbox{on $\RR^d.$} Then there exists $C:=c(\alpha)$ such that
\begin{eqnarray*}
\big{\|}|\textnormal{D}|^\alpha(f\circ\psi)-(|\textnormal{D}|^\alpha f)\circ\psi\big{\|}_{L^p}&\leq&
 C\max\Big(\vert 1-\|\nabla\psi^{-1}\|_{L^\infty}^{d+\alpha}\vert ; \vert 1-
 {\| \nabla\psi \|_{{L^\infty}}^{-d-\alpha}}\vert\Big) \\
 &&
 \|\nabla \psi \|_{L^\infty}^\alpha\| f \|_{\dot{B}_{p,1}^\alpha}.
\end{eqnarray*}
\end{prop}
\begin{proof} First we rule out the obvious  case $\alpha=0$ and let us recall
 the following formula detailed in \cite{C-C} which tells us that for all $\alpha\in]0,2[$
$$
|\textnormal{D}|^\alpha f(x)=C_{\alpha}\,\textnormal{P. V.}\int\frac{f(x)-f(y)}{\vert x-y \vert^{d+\alpha}}dy.
$$ 
Now we claim from (\ref{equivalence}) that  if $g\in\dot{B}_{p,1}^\alpha,$ with $\alpha\in]0,1[,$ then the above identity holds as an $L^p$ equality\begin{equation}\label{E1}
|\textnormal{D}|^\alpha f(x)=C_{\alpha}\int_{\RR^d}\frac{f(x)-f(y)}{\vert x-y \vert^{d+\alpha}}dy,\quad a.e.w.
\end{equation}
and moreover,
\begin{equation}\label{E2}
\||\textnormal{D}|^\alpha f\|_{L^p}\lesssim \|f\|_{\dot{B}_{p,1}^\alpha}.
\end{equation}
Indeed,  the $L^p$ norm of the integral function satisfies in view of Minkowski inequalities
$$
\big{\|}\int_{\RR^d}\frac{f(\cdot)-f(y)}{\vert \cdot-y \vert^{d+\alpha}}dy\big{\|}_{L^p}\leq
 \int_{\RR^d}\frac{\| f(\cdot)-f(\cdot-y)\|_{L^p}}{\vert y \vert^{d+\alpha}}dy\approx \|f\|_{\dot{B}_{p,1}^\alpha}.$$
 Thus we find that the left integral term is finite almost every where. \\
 Inasmuch as the flow preserves Lebesgue measure then the formula  (\ref{E1})  yields
  $$
 (|\textnormal{D}|^\alpha f)\circ\psi(x)=C_{\alpha}\int_{\RR^d}\frac{f(\psi(x))-f(y)}{\vert \psi(x)-y \vert^{d+\alpha}}dy=
 C_{\alpha}\int_{\RR^d}
 \frac{ f(\psi(x))-f(\psi(y))}{\vert \psi(x)-\psi(y) \vert^{d+\alpha}}dy 
 . $$
 Applying again (\ref{E1})
 with $f\circ\psi$, we obtain
 $$
 |\textnormal{D}|^\alpha (f\circ\psi)(x)=C_{\alpha}\int_{\RR^d}\frac{f(\psi(x))-f(\psi(y))}{\vert x-y \vert^{d+\alpha}}dy. 
 $$
 
 Thus we get
 \begin{eqnarray*}
 |\textnormal{D}|^\alpha(f\circ\psi)(x)-(|\textnormal{D}|^\alpha f)\circ\psi(x)
&=& C_{\alpha}\int_{\RR^d}\frac{f(\psi(x))-f(\psi(y))}{\vert x-y \vert^{d+\alpha}}\times\\
&&\Big( 1-\frac{\vert x-y \vert^{d+\alpha}}{\vert \psi(x)-\psi(y) \vert^{d+\alpha}} \Big) dy.
\end{eqnarray*}
 Taking the $L^p$ norm and using (\ref{equivalence}) we obtain
 \begin{equation}\label{hamouda1}
 \big{\|}|\textnormal{D}|^\alpha(f\circ\psi)-(|\textnormal{D}|^\alpha f)\circ\psi\big{\|}_{L^p}\lesssim 
 \| f \circ\psi\|_{\dot{B}_{p,1}^\alpha}
 \sup_{x,y}\big{\vert} 1-\frac{\vert x-y \vert^{d+\alpha}}{\vert \psi(x)-\psi(y) \vert^{d+\alpha}} \big{\vert}.
  \end{equation}
According to \cite{Marcel} one has  the following  composition result  
 $$
 \|f\circ\psi\|_{\dot{B}_{p,1}^\alpha}\leq c_{\alpha}\|\nabla\psi \|_{L^\infty}^\alpha \|f\|_{\dot{B}_{p,1}^\alpha}, \quad\textnormal{for}\quad\alpha\in]0,1[.
 $$
Therefore (\ref{hamouda1}) becomes
   \begin{eqnarray*}    \big{\|}|\textnormal{D}|^\alpha(f\circ\psi)-(|\textnormal{D}|^\alpha f)\circ\psi\big{\|}_{L^p}  & \leq &C\|\nabla \psi \|_{L^\infty}^\alpha
 \| f \|_{\dot{B}_{p,1}^\alpha} \times\\
 &&\sup_{x,y}\big{\vert} 1-\frac{\vert x-y \vert^{d+\alpha}}{\vert \psi(x)-\psi(y) \vert^{d+\alpha}}\big{\vert}.
 \end{eqnarray*}
It is plain from
 mean value Theorem that
 $$
\frac{1}{\|\nabla \psi\|_{L^\infty}^{d+\alpha}}\leq\frac{|x-y|^{d+\alpha}}{|\psi(x)-\psi(y)|^{d+\alpha}}\leq \|\nabla\psi^{-1}\|_{L^\infty}^{d+\alpha},
$$ 
 
which gives  easily the inequality
$$
 \sup_{x,y}\big{\vert} 1-\frac{\vert x-y \vert^{d+\alpha}}{\vert \psi(x)-\psi(y) \vert^{d+\alpha}} \big{\vert}
 \leq \max\Big(\vert 1-\|\nabla \psi^{-1} \|_{{L^\infty}}^{d+\alpha}\vert;
 \vert 1-{\|\nabla \psi \|_{{L^\infty}}^{-d-\alpha}}\vert \Big). 
  $$
  This concludes the proof.
  \end{proof}
  %%%%%%%%%%%%%%%%%%%%%%%%%%%%%%%%%%%%%%%%%
  %%%%%%%%%%%%%   Proof THeorem 1%%%%%%%%%
  %%%%%%%%%%%%
  %%%%%%%%%
   
  \section{Proof of Theorem \ref{Thm3}}
We shall divide our analysis into two cases: $r=+\infty$ and $r$ is finite. The first case is more easy and simply based upon
   a maximum principle and a commutator estimate. Before we move on let us mention  that in what follows we will work with the homogeneous Littlewood-Paley operators but we  take the same notation of the inhomogeneous operators.
   
    Set  $\theta_{q}:=\Delta_{q}\theta,$ then    localizing the  \mbox{$\QG$} equation through the operator $\Delta_{q}$ gives
 \begin{equation}\label{eq:1}
 \partial_{t}\theta_{q}+v\cdot\nabla\theta_{q}+|\hbox{D}|^{\alpha}\theta_{q}=-[\Delta_{q},v\cdot\nabla]\theta+f_{q}:=\mathcal{R}_{q}.
 \end{equation}
According to Proposition \ref{maximum} we have
\begin{equation}\label{eq:2}
\| \theta_{q}(t)\|_{L^p}\leq\|\theta_{q}^0\|_{L^p}+\int_{0}^t\|\mathcal{R}_{q}(\tau)\|_{L^p}d\tau.
\end{equation}
 Multiplying both sides by $2^{qs}$ and summing over $q$
  $$
 \|\theta\|_{\widetilde L^\infty_{t}\dot{B}_{p,1}^s}\leq \|\theta^0\|_{\dot{B}_{p,1}^s}+\|f\|_{L^1_{t}\dot{B}_{p,1}^s}+\int_{0}^t\sum_{q}2^{qs} \|\mathcal{R}_{q}(\tau) \|_{L^p}d\tau. $$
This yields in view of Lemma \ref{lemm12} 
 $$ 
 \|\theta\|_{\widetilde L^\infty_{t}\dot{B}_{p,1}^s}\leq \|\theta^0\|_{\dot{B}_{p,1}^s}+\|f\|_{L^1_{t}\dot{B}_{p,1}^s}+C\int_{0}^t\|\nabla v(\tau)\|_{L^\infty} \|\theta\|_{\widetilde L^\infty_{\tau}\dot{B}_{p,1}^s}d\tau. $$
To achieve the proof in the case of $r=\infty,$  it suffices to use Gronwall's inequality.\\

We shall now turn to the proof of the finite case $r<\infty$ which is more technical.
 Let  $\psi$ denote the flow of the  velocity $v$ and  set 
 $$
 \bar\theta_{q}(t,x)=\theta_{q}(t,\psi(t,x))\quad\hbox{and}\quad \bar {\mathcal R}_{q}(t,x)=\mathcal{R}_{q}(t,\psi(t,x)).
 $$

 Since the flow preserves Lebesgue measure then we obtain
 \begin{equation}
 \label{eq:5}
\|\bar{\mathcal{R}}_{q}\|_{L^p}\leq \|[\Delta_{q},v\cdot\nabla]\theta\|_{L^p}+\|f_{q}\|_{L^p}.
\end{equation}
It is not hard to check that the function $\bar\theta_{q}$ satisfies 
 \begin{equation}
 \label{T1}
 \partial_{t}\bar\theta_{q}+|\hbox{D}|^\alpha\bar\theta_{q}=
 |\hbox{D}|^\alpha(\theta_{q}\circ\psi)-(|\hbox{D}|^\alpha\theta_{q})\circ\psi+\bar {\mathcal R}_{q}:=\bar{\mathcal R}_{q}^1.
 \end{equation}
From Proposition \ref{pr:1} we find that for  $q\in\ZZ$
\begin{eqnarray}
\label{T:5}
\nonumber\| |\hbox{D}|^\alpha(\theta_{q}\circ\psi)-(|\hbox{D}|^\alpha\theta_{q})\circ\psi\|_{L^p}
 &\leq& Ce^{CV(t)} \Big(e^{CV(t)}-1\Big)\|\theta_{q}(t)\|_{\dot B_{p,1}^\alpha}\\
 &\leq&Ce^{CV(t)} \Big(e^{CV(t)}-1\Big) 2^{q\alpha}\|\theta_{q}\|_{L^p},
  \end{eqnarray}
where $V(t):=\|\nabla v\|_{L^1_{t}L^\infty}.$ Notice that we have used here the classical estimates
 $$
 e^{-CV(t)}\leq\|\nabla\psi^{\mp 1}(t)\|_{L^\infty}\leq
  e^{CV(t)}.
   $$
Putting together (\ref{eq:5}) and (\ref{T:5}) yield
$$
\|\bar{\mathcal R}_{q}^1(t)\|_{L^p}\leq \| f_{q}(t)\|_{L^p}+\|[\Delta_{q},v\cdot\nabla]\theta\|_{L^p}
+Ce^{CV(t)
} (e^{CV(t)}-1)2^{q\alpha}\|\theta_{q}(t)\|_{L^p}.$$
  Applying the operator $\Delta_{j},$ for $j\in\ZZ,$ to the equation (\ref{T1}) and using Proposition \ref{l:5}
\begin{eqnarray} 
 \|\Delta_{j}\bar\theta_{q}(t)\|_{L^p}&\leq& Ce^{-ct2^{j\alpha}}\|\Delta_{j}\theta_{q}^0\|_{L^p}+C\int_{0}^te^{-c (t-\tau)2^{j\alpha}}\| f_{q}(\tau)\|_{L^p}d\tau \\
\nonumber&+& Ce^{CV(t)}(e^{CV(t)}-1)2^{q\alpha}
 \int_{0}^te^{-c (t-\tau)2^{j\alpha}}\|\theta_{q}(\tau)\|_{L^p}d\tau\\
\nonumber &+&C\int_{0}^te^{-c (t-\tau)2^{j\alpha}}\|[\Delta_{q},v\cdot\nabla]\theta(\tau)\|_{L^p}d\tau.
  \end{eqnarray}
  Integrating this estimate with respect to the time  and using Young's inequality 
\begin{eqnarray}\label{eres}
\nonumber \|\Delta_{j}\bar\theta_{q}\|_{L^r_{t}L^p}&\leq& C
2^{-j\alpha/r}\big((1-e^{-crt2^{j\alpha}})^{\frac{1}{r}}\|\Delta_{j}\theta_{q}^0\|_{L^p}+\|f_{q}\|_{L^1_{t}L^p}\big)\\
\nonumber&+&
Ce^{CV(t)}(e^{CV(t)}-1)2^{(q-j)\alpha}\|\theta_{q}\|_{L^r_{t}L^p}\\
&+&C2^{-j\alpha/r}\int_{0}^t \|[\Delta_{q},v\cdot\nabla]\theta(\tau)\|_{L^p}d\tau.
\end{eqnarray}
Since the flow $\psi$ preserves  Lebesgue measure then one writes
\begin{eqnarray*}2^{q(s+\alpha/r)}\|\theta_{q}\|_{L^r_{t}L^p}&=&2^{q(s+\alpha/r)}
\|\bar\theta_{q}\|_{L^r_{t}L^p}\\
&\leq &2^{q(s+\alpha/r)}\Big(
\sum_{\vert j-q \vert >N}\|\Delta_{j}\bar\theta_{q}\|_{L^r_{t}L^p}+
\sum_{\vert j-q \vert\leq N}\|\Delta_{j}\bar\theta_{q}\|_{L^r_{t}L^p}\Big)\\
&:=&\textnormal{I}+\textnormal{II}.
\end{eqnarray*}
To estimate the term $\textnormal I$ we make appeal to  Lemma \ref{l400} 
\begin{eqnarray*}
\|\Delta_{j}\bar\theta_{q}\|_{L^r_{t}L^p}&\leq& C2^{-\vert q-j\vert}
e^{\int_{0}^t
\|\nabla v(\tau)\|_{L^\infty}d\tau}\|\theta_{q}\|_{L^r_{t}L^p}\\
&\leq&C2^{-\vert q-j\vert}
e^{V(t)}\|\theta_{q}\|_{L^r_{t}L^p}.
\end{eqnarray*}
Therefore we get
\begin{equation}\label{e:11}
\textnormal{I}\leq C2^{-N}e^{V(t)}2^{q(s+\alpha/r)}\|\theta_{q}\|_{L^r_{t}L^p}.
\end{equation}
In order to bound the second term $\textnormal{II}$ we use (\ref{eres}) 
\begin{eqnarray}\label{Meth1}
\nonumber \textnormal{II}&\leq& C(1-e^{-crt2^{q\alpha}})^{\frac{1}{r}}2^{qs}\|\theta_{q}^0\|_{L^p}+C2^{N\frac{\alpha}{r}}2^{qs}\|f_{q}\|_{L^1_{t}L^p}\\
\nonumber&+&
C2^{N\alpha}e^{CV(t)}(e^{CV(t)}-1)2^{q(s+\alpha/r)}\|\theta_{q}\|_{L^r_{t}L^p}\\
&+&C2^{N\alpha/r}2^{qs}\int_{0}^t \|[\Delta_{q},v\cdot\nabla]\theta(\tau)\|_{L^p}d\tau.
\end{eqnarray}
Denote  $Z_{q}^r(t):=2^{q(s+\alpha/r)}\|\theta_{q}\|_{L^r_{t}L^p},$ then we obtain  in view of  (\ref{e:11}) and (\ref{Meth1}) 
\begin{eqnarray*}
Z_{q}^r(t)&\leq &C(1-e^{-crt2^{q\alpha}})^{\frac{1}{r}}2^{qs}\|\theta_{q}^0\|_{L^p}+C2^{N\frac{\alpha}{r}}2^{qs}\|f_{q}\|_{L^1_{t}L^p}\\
&+&C\big(2^{N\alpha}e^{CV(t)}
(e^{CV(t)}-1)+2^{-N}e^{CV(t)}\big)Z_{q}^r(t)\\
&+&C2^{N\alpha/r}2^{qs}\int_{0}^t \|[\Delta_{q},v\cdot\nabla]\theta(\tau)\|_{L^p}d\tau.
\end{eqnarray*}
We can easily show that there exists two pure constants $N$ and $C_{0}$ such that
$$
V(t)\leq C_{0}\Rightarrow C2^{-N}e^{CV(t)}+C2^{N\alpha}e^{CV(t)}(e^{CV(t)}-1)\leq \frac{1}{2}\cdot
$$ 
Thus we obtain under this condition
\begin{eqnarray}
\label{mahma00}
\nonumber Z_{q}^r(t)&\leq& C(1-e^{-crt2^{q\alpha}})^{\frac{1}{r}}2^{qs}\|\theta_{q}^0\|_{L^p}+C2^{qs}\|f_{q}\|_{L^1_{t}L^p}\\
&+&C2^{qs}\int_{0}^t \|[\Delta_{q},v\cdot\nabla]\theta(\tau)\|_{L^p}d\tau.
\end{eqnarray} 
Summing over $q$  and using Lemma \ref{lemm12} lead  for $V(t)\leq C_{0},$
\begin{eqnarray*}
 \|\theta\|_{\widetilde L^r_{t}\dot{B}_{p,1}^{s+\frac{\alpha}{r}}}&\leq& C \|\theta^0\|_{\dot{B}_{p,1}^s}+C\|f\|_{L^1_{t}\dot B_{p,1}^s}+
 C\int_{0}^t\| \nabla v(\tau)\|_{L^\infty}
\|\theta(\tau)\|_{\dot{B}_{p,1}^s}d\tau\\
&\leq&C \|\theta^0\|_{\dot{B}_{p,1}^s}+C\|f\|_{L^1_{t}\dot B_{p,1}^s}+CV(t)\| \theta\|_{L^\infty_{t}\dot{B}_{p,1}^s}.
\end{eqnarray*}
Thus we get in view of the estimate of the case $r=\infty$
\begin{equation}\label{EDF}
 \|\theta\|_{\widetilde L^r_{t}\dot{B}_{p,1}^{s+\frac{\alpha}{r}}}\leq C \|\theta^0\|_{\dot{B}_{p,1}^s}+C\|f\|_{L^1_{t}\dot B_{p,1}^s}.
 \end{equation}
This gives the result for a short time. 

For an arbitrary positive time  $T$ we make a partition $(T_{i})_{i=0}^M$  of the interval $[0,T],$ 
such \mbox{that $\displaystyle \int_{T_{i}}^{T_{i+1}}\|\nabla v(\tau)\|_{L^\infty}d\tau\approx C_{0}.$}
Then proceeding for (\ref{EDF}), we obtain 
$$ \|\theta\|_{\widetilde L^r_{[T_{i},T_{i+1}]}\dot{B}_{p,1}^{s+\frac{\alpha}{r}}}\leq
 C \|\theta(T_{i})\|_{\dot{B}_{p,1}^s}+C\int_{T_{i}}^{T_{i+1}}\|f(\tau)\|_{\dot{B}_{p,1}^s}d\tau.$$
Applying the triangle inequality gives
$$
\|\theta\|_{\widetilde L^r_{T}\dot{B}_{p,1}^{s+\frac{\alpha}{r}}}\leq C\sum_{i=0}^{M-1}\|\theta(T_{i})\|_{\dot{B}_{p,1}^s}+C\int_{0}^T\|f(\tau)\|_{\dot B_{p,1}^s}d\tau.$$
On the other hand the estimate proven in the case $r=\infty$ allows us to write
$$
\|\theta\|_{\widetilde L^r_{T}\dot{B}_{p,1}^{s+\frac{\alpha}{r}}}\leq 
CM\big(\|\theta^0\|_{\dot{B}_{p,1}^s}+\|f\|_{L^1_{T}\dot B_{p,1}^s}\big)e^{CV(T)}+\|f\|_{L^1_{T}\dot B_{p,1}^s}.$$
Thus the following observation $C_{0}M\approx 1+V(t)$ completes the proof of the theorem.

 \section{Proof of Theorem {\ref{Thm1}}}
  For the sake of a concise presentation, we shall just provide the {\it a priori} estimates supporting the claims of the theorem. To achieve the proof one must combine in a standard way  these estimates with a standard
  approximation procedure such as the following iterative scheme
  $$
\left\lbrace
\begin{array}{l}\partial_{t}\theta_{n+1}+v^n\cdot\nabla\theta_{n+1}+|\hbox{D}|^{{\alpha}} \theta_{n+1}=0,\\
v_{n}=(-R_{2} \theta_{n},R_{1}\theta_{n}),\\
\theta_{n+1}(0,x)=S_{n}\theta^0(x),\\
(\theta_{0},v_{0})=(0,0).
\end{array}
\right.$$
%%%%%%
%%%%%% GLobal existence
%%%%%%
  \subsection{ Global  existence} It is plain from  Theorem \ref{Thm3} that to derive global {\it a priori} estimates it is sufficient to bound globally in time  the quantity 
 $
 V(t):=\|\nabla v\|_{L^1_{t}L^\infty}. $  First, 
 the embedding $\dot{B}_{\infty,1}^0\hookrightarrow L^\infty$ combined with  the fact that Riesz transform maps continuously  homogeneous Besov space into itself 
\begin{equation}\label{important}
\|\nabla v\|_{L^1_{t}L^\infty}\leq \|\nabla v\|_{L^1_{t}\dot B_{\infty,1}^0}\leq C\|\theta\|_{L^1_{t}\dot B_{\infty,1}^1}.
\end{equation}
Combined with Theorem \ref{Thm3} this yields
$$
V(t)\leq C\|\theta^0\|_{\dot{B}_{\infty,1}^{1-\alpha}} e^{CV(t)}.
$$
 Since the function $V$ depends continuously in time and $V(0)=0$ then we can deduce that for small initial data $V$ does not blow up, and there exists $C_{1},\eta>0$ such that
   \begin{equation}\label{TA1}
 \|\theta^0\|_{\dot B_{\infty,1}^{1-\alpha}}<\eta\Rightarrow\| \nabla v\|_{L^1(\RR_{+};L^\infty)}\leq C_{1}\|\theta^0\|_{\dot{B}_{\infty,1}^{1-\alpha}},\forall t\in\RR_{+}.  \end{equation}
  Let us now show how to derive the {\it a priori} estimates. Take  $s\geq s_{c}^p:=1+\frac{2}{p}-\alpha.$ Then combining Theorem \ref{Thm3} with (\ref{TA1}) we get
  \begin{eqnarray*}
   \|\theta\|_{\widetilde L^\infty_{\RR_{+}}\dot{B}_{p,1}^{s}}+ \|\theta\|_{ L^1_{\RR_{+}}\dot{B}_{p,1}^{s+\alpha}}    &\leq& C\|\theta^0\|_{\dot{B}_{p,1}^s}e^{C\|\theta^0\|_{\dot{B}_{\infty,1}^{1-{\alpha}}} }\\
    &\leq&
    C\|\theta^0\|_{\dot{B}_{p,1}^s}. 
        \end{eqnarray*}
    On the other hand we have from Proposition \ref{maximum}
    $$
 \forall t\in\RR_{+},\,   \|\theta(t)\|_{L^p}\leq\|\theta^0\|_{L^p}. 
    $$
    Therefore we get an estimate of $\theta$ in the  inhomogeneous Besov space as follows
    $$
\|\theta\|_{\widetilde L^\infty_{\RR_{+}}{B}_{p,1}^{s}}\leq C\|\theta^0\|_{{B}_{p,1}^s}.    
$$
 Using again Theorem \ref{Thm3} yields
 $$
 \|\theta\|_{\widetilde L^\infty_{t}\dot{B}_{\infty,1}^0}\leq C\|\theta^0\|_{\dot{B}_{\infty,1}^0}e^{CV(t)}\leq
 C\|\theta^0\|_{\dot{B}_{\infty,1}^0},
 $$
Thus we obtain for $p\in[1,\infty]$
\begin{equation}\label{Jui1}
\|\theta\|_{\widetilde L^\infty_{\RR_{+}}{\mathcal{X}}_{p}^s}\leq C\|\theta^0\|_{{\mathcal{X}}_{p}^s}.
\end{equation}    For the velocity we have the following result.
\begin{lem}
For $p\in]1,\infty]$ there exists $C_{p}$ such that
$$
\|v\|_{\widetilde L^\infty_{\RR_+}B_{p,1}^s}\leq C_{p}\|\theta^0\|_{\mathcal{X}_{p}^s}.
$$
However, for $p=1$ we have
$$
\|v\|_{\widetilde L^\infty_{\RR_{+}}\dot B_{1,1}^s}+\|v\|_{L^\infty_{\RR_{+}}L^{p_{1}}}\leq C_{p_{1}}\|\theta^0\|_{B_{1,1}^s},\,\forall p_{1}>1.
$$
\end{lem}
\begin{proof}
Let $p\in]1,\infty[.$ Then we can write in view of (\ref{Jui1})
     \begin{eqnarray*}
     \|v\|_{\widetilde L^\infty_{\RR_{+}}B_{p,1}^s}&\leq &\|v\|_{\widetilde L^\infty_{\RR_{+}}\dot{B}_{p,1}^s}+\|\Delta_{-1}v
\|_{L^\infty_{\RR_{+}}L^p}\\
&\leq& C\|\theta\|_{\widetilde L^\infty_{\RR_{+}}\dot{B}_{p,1}^s}+C\|v\|_{L^\infty_{\RR_{+}}L^p}\\
&\leq&C\|\theta^0\|_{B_{p,1}^s}+C\|v\|_{L^\infty_{\RR_{+}}L^p}.
\end{eqnarray*}
 Combining the boundedness of Riesz transform  with the maximum principle
$$
\|v\|_{L^\infty_{\RR_{+}}L^p}\leq C_{p}\|\theta^0\|_{L^p}.
$$
Thus we obtain
$$
\|v\|_{\widetilde L^\infty_{\RR_{+}}B_{p,1}^s}\leq C_{p}\|\theta^0\|_{B_{p,1}^s}.
$$
To treat  the case  $p=\infty$ we write  according to the embedding $\dot{B}_{\infty,1}^0\hookrightarrow  L^\infty$ and the continuity of Riesz transform
 $$
 \|\Delta_{-1}v(t)\|_{L^\infty}\leq C\|\theta(t)\|_{\dot{B}_{\infty,1}^0}. $$
Combining this estimate with (\ref{Jui1}) yields
 $$
 \|v\|_{\widetilde L^\infty_{\RR_{+}}B_{\infty,1}^s}\leq C\|\theta^0\|_{B_{\infty,1}^s\cap \dot{B}_{\infty,1}^0}.
 $$
Hence we get for all $p\in]1,\infty]$
\begin{equation}\label{acigne}
\|v\|_{\widetilde L^\infty_{\RR_{+}}B_{p,1}^s}\leq C_{p}\|\theta^0\|_{{\mathcal{X}}_{p}^s}.
\end{equation}
Let us now move to the case $p=1.$ Since $B_{1,1}^s\hookrightarrow  L^{p_{1}}$ for all $p_{1}\geq1$ then we get in view of Bernstein's inequality and the maximum principle
$$
\|\Delta_{-1}v\|_{L^{p_{1}}}\leq C_{p_{1}}\|\theta^0\|_{L^{p_{1}}}\leq C_{p_{1}}\|\theta^0\|_{B_{1,1}^s}.
$$
We eventually find that $v\in\widetilde L^\infty_{\RR_{+}}\dot B_{1,1}^s\cap L^\infty_{\RR_{+}} L^{p_{1}}.$

\end{proof}
Let us now briefly sketch the proof of the continuity in time, that is      $\theta\in C(\RR_{+};{\mathcal{X}}_{p}^s).$ We should  only treat the  finite case of  $p$ and similarly one can show the case $p=\infty$. From the definition of Besov spaces we have
$$
\|\theta(t)-\theta(t')\|_{B_{p,1}^s}\leq \sum_{q< N}2^{qs}\|\theta_{q}(t)-\theta_{q}(t')\|_{L^p}+2\sum_{q\geq N}2^{qs}\|\theta_{q}\|_{L^\infty_{\RR_{+}}L^p}$$
Let $\epsilon>0$ then we get from  (\ref{Jui1}) the existence of a number $N$ such that
$$
\sum_{q\geq N}2^{qs}\|\theta_{q}\|_{L^\infty_{\RR_{+}}L^p}\leq\frac{\epsilon}{4}.
$$
Thanks to Taylor's formula
\begin{eqnarray*}
\sum_{q< N}2^{qs}\|\theta_{q}(t)-\theta_{q}(t')\|_{L^p}&\leq &|t-t'|\sum_{q<N}2^{qs}\|\partial_{t}\theta_{q}\|_{L^\infty_{\RR_{+}}L^p}\\
&\leq& C|t-t'|2^N\|\partial_{t}\theta\|_{L^\infty_{\RR_{+}}B_{p,1}^{s-1}}.
\end{eqnarray*}
To estimate the last term
 we write
$$
\partial_{t}\theta=-|\textnormal{D}|^\alpha\theta-v\cdot\nabla\theta.
$$
In one hand we have $|\textnormal{D}|^\alpha\theta\in B_{p,1}^{s-\alpha}\hookrightarrow B_{p,1}^{s-1}$. On the other hand since the space $B_{p,1}^s$ is an algebra ($s>\frac{2}{p}$) and $v$ is zero divergence
then
 $$
\|v\cdot\nabla\theta\|_{B_{p,1}^{s-1}}\leq C\|v\,\theta\|_{B_{p,1}^{s}}\leq C\|v\|_{B_{p,1}^{s}}\|\theta\|_{B_{p,1}^{s}}.
$$
Thus we get $\partial_{t}\theta\in L^\infty_{\RR_{+}}B_{p,1}^{s-1}$ and this allows us to finish the proof of the continuity.

%%%%%%
%%%% LOCAL ExiSTENCE
%%%%%%%
%%%%%%%
%%%%%%%%%
 \subsection { Local existence} The local time existence depends on the control of the quantity $V(t):=\|\nabla v\|_{L^1_{t}L^\infty}.$ In our analysis we distinguish two cases:

 $\bullet$ {\it First case: } $s>s_{c}^p=1+\frac{2}{p}-\alpha.$ 
 
 We observe first that there exists $r>1$ such that $1+\frac{2}{p}-\frac{\alpha}{r}\leq s.$
 From  (\ref{important}) and according to the said H\"older's inequality we have
  \begin{eqnarray*}
  V(t)&\leq& C\|\theta\|_{L^1_{t}\dot{B}_{\infty,1}^1}\\
  &\leq&Ct^{\frac{1}{\bar{r}}} \|\theta\|_{L^r_{t}\dot{B}_{\infty,1}^1}.
    \end{eqnarray*}
Using Theorem (\ref{Thm3}) we obtain 
 \begin{eqnarray*}
V(t)\leq Ct^{\frac{1}{\bar{r}}} 
 \|\theta^0\|_{\dot{B}_{\infty,1}^{1-\frac{\alpha}{r}}}e^{CV(t)}.
 \end{eqnarray*}
  Thus we conclude that there exists $C_{0},\,\eta>0$  such that 
  \begin{equation}\label{Tr23}
  t^{\frac{1}{\bar{r}}} 
 \|\theta^0\|_{\dot{B}_{\infty,1}^{1-\frac{\alpha}{r}}}\leq\eta\Rightarrow V(t)\leq C_{0},
   \end{equation}
 and this
   gives from Theorem \ref{Thm3}
\begin{equation}\label{Tr24}
 \|\theta\|_{L^\infty_{t
 }B_{p,1}^s}+\|\theta\|_{L^1_{t}\dot B_{p,1}^{1+\frac{2}{p}}}\leq C\|\theta^0\|_{B_{p,1}^s}.
 \end{equation}
   We point out that one can deduce from (\ref{Tr23}) that the time  existence is bounded below  
   $$T\gtrsim \|\theta^0\|_{\dot{B}_{\infty,1}^{1-\frac{\alpha}{ r}}}^{-\bar{r}}.$$  
    
     $\bullet$ {\it Second case: }$ s=s_{c}^p=1+\frac{2}{p}-\alpha.$
     
  By applying (\ref{mahma00}) to the $\QG$ equation with $r=1,\,p=\infty$ and $s=1-\alpha$ we have under the condition $V(t)\leq C_{0}$
  $$ 
  \|\theta\|_{L^1_{t}\dot{B}_{\infty,1}^1}\leq C\sum_{q\in\ZZ}(1-e^{-ct2^{q\alpha}})2^{q(1-\alpha)}\|\theta_{q}^0\|_{L^\infty}+C\sum_{q\in\ZZ}2^{q(1-\alpha)}\|[\Delta_{q},v\cdot\nabla]\theta\|_{L^1_{t}L^\infty}.
  $$   
The second term of the right-hand side can be estimated from  Lemma \ref{lemm12} as follows
\begin{eqnarray}
\label{All-Port}
\nonumber\sum_{q\in\ZZ}2^{q(1-\alpha)}\|[\Delta_{q},v\cdot\nabla]\theta\|_{L^1_{t}L^\infty}&\leq& C\|v\|_{\widetilde L^2_{t}\dot B_{\infty,1}^{1-\frac{\alpha}{2}}}\|\theta\|_{\widetilde L^2_{t}\dot B_{\infty,1}^{1-\frac{\alpha}{2}}}\\
&\leq&
C\|\theta\|_{\widetilde L^2_{t}\dot B_{\infty,1}^{1-\frac{\alpha}{2}}}^2.
\end{eqnarray}
Notice that we have used in the above inequality the fact that Riesz transform maps continuously  homogeneous Besov space into itself. Hence we get
\begin{equation}
\label{RT}
 \|\theta\|_{L^1_{t}\dot{B}_{\infty,1}^1}\leq C\sum_{q\in\ZZ}(1-e^{-ct2^{q\alpha}})2^{q(1-\alpha)}\|\theta_{q}^0\|_{L^\infty}+C\|\theta\|_{\widetilde L^2_{t}\dot B_{\infty,1}^{1-\frac{\alpha}{2}}}^2.
 \end{equation}
Using again (\ref{mahma00}) with $r=2,\,p=\infty$ and $s=1-{\alpha},$ we obtain
$$
\|\theta\|_{\widetilde L^2_{t}\dot B_{\infty,1}^{1-\frac{\alpha}{2}}}\leq
C\sum_{q\in\ZZ}(1-e^{-ct2^{q\alpha}})^{\frac{1}{2}}2^{q(1-\alpha)}\|\theta_{q}^0\|_{L^\infty}+C\sum_{q\in\ZZ}2^{q(1-\alpha)}\|[\Delta_{q},v\cdot\nabla]\theta\|_{L^1_{t}L^\infty}.
$$
 Thus (\ref{All-Port}) yields
 $$
 \|\theta\|_{\widetilde L^2_{t}\dot B_{\infty,1}^{1-\frac{\alpha}{2}}}\leq
C\sum_{q\in\ZZ}(1-e^{-ct2^{q\alpha}})^{\frac{1}{2}}2^{q(1-\alpha)}\|\theta_{q}^0\|_{L^\infty}+C \|\theta\|_{\widetilde L^2_{t}\dot B_{\infty,1}^{1-\frac{\alpha}{2}}}^2.
 $$   
By Lebesgue theorem we have 
$$
\lim_{t\to 0^+}\sum_{q\in\ZZ}(1-e^{-ct2^{q\alpha}})^{\frac{1}{2}}2^{q(1-\alpha)}\|\theta_{q}^0\|_{L^\infty}=0.$$
Let $\eta$ be a sufficiently small constant and define
$$
T_{0}:=\sup\Big\{t> 0,\sum_{q\in\ZZ}(1-e^{-ct2^{q\alpha}})^{\frac{1}{2}}2^{q(1-\alpha)}\|\theta_{q}^0\|_{L^\infty}\leq \eta\Big\}.
$$
 Then we have under the assumptions  $t\leq T_{0}$ and  $V(t)\leq C_{0}$ 
 $$
 \|\theta\|_{\widetilde L^2_{t}\dot B_{\infty,1}^{1-\frac{\alpha}{2}}}\leq 2C\sum_{q\in\ZZ}(1-e^{-ct2^{q\alpha}})^{\frac{1}{2}}2^{q(1-\alpha)}\|\theta_{q}^0\|_{L^\infty}.  $$
 Inserting this estimate into (\ref{RT}) gives 
 \begin{eqnarray*}
 V(t)\leq C\|\theta\|_{L^1_{t}\dot{B}_{\infty,1}^1}&\leq&
 C\sum_{q\in\ZZ}(1-e^{-ct2^{q\alpha}})^{\frac{1}{2}}2^{q(1-\alpha)}\|\theta_{q}^0\|_{L^\infty}\\
 &+&C\Big(\sum_{q\in\ZZ}(1-e^{-ct2^{q\alpha}})^{\frac{1}{2}}2^{q(1-\alpha)}\|\theta_{q}^0\|_{L^\infty}\Big)^2. 
 \end{eqnarray*}
 For sufficiently small $\eta$ we obtain $V(t)<C_{0}$ and this allows us to prove  that the time $T_{0}$ is actually  a local time existence.
 Thus we obtain from Theorem \ref{Thm3}
 $$
\|\theta\|_{\widetilde L^\infty_{T}B_{p,1}^{s_{c}^p}}+\|\theta\|_{L^1_{T}\dot{B}_{p,1}^{1+\frac{2}{p}}}\leq C\|\theta^0\|_{B_{p,1}^{s_{c}^p}}.
 $$
  %%
 %%%%%% UNIQUENSSSSSSSSSS
 %%%%%%%%
 %%%%%%%%%  
   %%%%%%%%%
   %%%%%
   %%%%
    \subsection{ Uniqueness}
  We shall give the proof of the uniqueness result which can be formulated as follows. There exists at most one solution for the system $\QG$ in the functions space $X_{T}:=L^\infty_{T}\dot{B}_{\infty,1}^0\cap L^1_{T}\dot{B}_{\infty,1}^1.$ We stress out that the space $L^\infty_{T}X^s_{p}\cap L^1_{T}\dot B_{p,1}^{s+\alpha},$ with $p\in[1,\infty],$ is continuously embedded in $X_{T}.$\\
   Let $\theta^{i}, i=1,2$ (and $v^{i}$ the corresponding velocity) be two solutions of the $\QG$ equation with the same  initial data  and  belonging to the space $X_{T}.$ We set $\theta=\theta^{1}-\theta^{2}$ and $v=v^1-v^2,$ then it is plain that
  $$
\partial_{t}\theta+v^1\cdot\nabla\theta+|\textnormal{D}|^\alpha \theta=-v\cdot\nabla\theta^{2},\,\,\theta_{|t=0}=0.
$$  
Applying Theorem \ref{Thm3} to this equation gives
 \begin{equation}\label{Coupe1}
 \|\theta(t)\|_{\dot{B}_{\infty,1}^0}\leq Ce^{C\|\nabla v^1\|_{L^1_{t}L^\infty}} \int_{0}^t\|v\cdot\nabla\theta^2(\tau)\|_{\dot{B}_{\infty,1}^0}d\tau.
 \end{equation}
We will now make use of the following law product and its proof will be given later.
\begin{equation}\label{Ac1}
\|v\cdot\nabla\theta^2\|_{\dot{B}_{\infty,1}^0}\leq C \|v\|_{\dot{B}_{\infty,1}^0}\|\theta^2\|_{\dot{B}_{\infty,1}^1}.
\end{equation}
Since Riesz transform maps continuously $\dot{B}_{\infty,1}^0$ into itself, then we get
$$
\|v\cdot\nabla\theta^2\|_{\dot{B}_{\infty,1}^0}\leq C \|\theta\|_{\dot{B}_{\infty,1}^0}\|\theta^2\|_{\dot{B}_{\infty,1}^1}.$$   
 Inserting this estimate into (\ref{Coupe1}) and using  Gronwall's inequality give  the wanted result.\\
 Let us now turn to the proof of (\ref{Ac1}) which is based on Bony's decomposition
 $$
 v\cdot\nabla\theta^2=T_{v}\nabla\theta^2+T_{\nabla\theta^2}v+R(v,\nabla\theta^2), \quad\hbox{with}
 $$
 $$
 T_{v}\nabla\theta^2=\sum_{q\in\ZZ}\dot S_{q-1}v\cdot\nabla\dot\Delta_{q}\theta^2\quad\hbox{and}\quad
 R(v,\nabla\theta^2)=\sum_{q\in\ZZ\atop
i\in\{\mp1,0\} }\dot\Delta_{q}v\cdot\dot\Delta_{q+i}\nabla\theta^2.
  $$ 
  Using the quasi-orthogonality of the paraproduct terms one obtains
  \begin{eqnarray*}
  \|T_{v}\nabla\theta^2\|_{\dot B_{\infty,1}^0}&\leq& C\sum_{q\in\ZZ}\|\dot S_{q-1}v\|_{L^\infty}\|\dot\Delta_{q}\nabla\theta^2\|_{L^\infty}\\
  &\leq&C\|v\|_{\dot B_{\infty,1}^0}\|\theta^2\|_{\dot B_{\infty,1}^1}.
  \end{eqnarray*}
  By the same way we get
   \begin{eqnarray*}
  \|T_{\nabla\theta^2}v\|_{\dot B_{\infty,1}^0}&\leq& C\sum_{q\in\ZZ}\|\dot S_{q-1}\nabla\theta^2\|_{L^\infty}\|\dot\Delta_{q}v\|_{L^\infty}\\
  &\leq& C\|\nabla\theta^2\|_{L^\infty}\|v\|_{\dot B_{\infty,1}^0}\\
  &\leq&C\|\theta^2\|_{\dot B_{\infty,1}^1}\|v\|_{\dot B_{\infty,1}^0}.
  \end{eqnarray*}
  For the remainder term we write in view of the incompressibility of the velocity and the convolution inequality
    \begin{eqnarray*}
  \|R(v,{\nabla\theta^2})\|_{\dot B_{\infty,1}^0}&=&\sum_{j\in\ZZ} \|\dot\Delta_{j}R(v,\nabla\theta^2)\|_{L^\infty} \leq C\sum_{q\geq j-3
  \atop i\in\{\mp1,0\}}2^j\|\dot \Delta_{q}v\|_{L^\infty}\|\dot\Delta_{q+i}\theta^2\|_{L^\infty}\\
  &\leq&C\sum_{q\geq j-3
  \atop i\in\{\mp 1,0 \}}2^{j-q}\|\dot \Delta_{q}v\|_{L^\infty}2^q\|\dot\Delta_{q+i}\theta^2\|_{L^\infty}\\
    &\leq& C\|\nabla\theta^2\|_{L^\infty}\|v\|_{\dot B_{\infty,1}^0}\\
  &\leq&C\|\theta^2\|_{\dot B_{\infty,1}^1}\|v\|_{\dot B_{\infty,1}^0}.
  \end{eqnarray*}
  This completes the proof of (\ref{Ac1}). 
   \section{Appendix}
   The following result is due to Vishik \cite{v1} and was used in a crucial way for the proof of Theorem {\ref{Thm3}}. For the convenience of the reader we will give a short  proof based on the duality method.
 \begin{lem}
 \label{l400} 
Let  $f$ be a function in  Schwartz class and $\psi$ a 
diffeomorphism   preserving  Lebesgue measure, then we have for all
$p\in[1,+\infty]$ and for all  { $j,q\in\ZZ,$}
$$\|\dot\Delta_j(\dot\Delta_q f\circ\psi)\|_{L^p}\leq
C2^{-\vert j-q\vert}\|\nabla\psi ^{\epsilon(j,q)}\|_{L^{\infty}}\|\dot\Delta_q
f\|_{L^p},$$
with
$$\epsilon(j,q)=\hbox{sign}(j-q).$$
\end{lem}
  We shall begin with the proof of Lemma \ref{l400}
   \begin{proof}
We distinguish two cases: $j\geq q$ and $j<q.$ For the first one we simply use Bernstein's inequality 
$$
\|\dot\Delta_j(\dot\Delta_q f\circ\psi)\|_{L^p}\lesssim 2^{-j}\|\nabla\dot\Delta_j(\dot\Delta_q f\circ\psi)\|_{L^p}.
$$
It suffices now to combine Leibnitz formula  again with Bernstein's inequality
\begin{eqnarray*}
\|\nabla\dot\Delta_j(\dot\Delta_q f\circ\psi)\|_{L^p}&\lesssim& \|\nabla\dot\Delta_q f\|_{L^p}\|\nabla\psi\|_{L^\infty}\\
&\lesssim& 2^q\|\dot\Delta_q f\|_{L^p}\|\nabla\psi\|_{L^\infty}.
\end{eqnarray*}
This yields to the desired inequality. Let us now move to the second case and use the following duality result
\begin{equation}
\label{duality1}
\|\dot\Delta_j(\dot\Delta_q f\circ\psi)\|_{L^p}=\sup_{\|g\|_{L^{\bar p}}\leq1}\big{|}\langle\dot\Delta_j(\dot\Delta_q f\circ\psi),g\rangle\big{|},\,\hbox{with}\quad \frac{1}{p}+\frac{1}{\bar p}=1.
\end{equation}
Let $\bar\varphi\in C^\infty_{0}(\RR^d)$ be supported in a ring and taking value $1$ on the ring $\mathcal{C}$ (see the definition of the dyadic decomposition). We set $\bar{\dot\Delta}_{q}f:=\bar\varphi(2^{-q}\textnormal{D})f.$ Then we can see easily that $\dot\Delta_{q}f=\bar{\dot\Delta}_{q}\dot\Delta_{q}f.$
Combining this fact with Parseval's identity and the preserving measure by the flow
$$
\big{|}\langle\dot\Delta_j(\dot\Delta_q f\circ\psi),g\rangle\big{|}=\big{|}\langle\dot\Delta_q f,\bar{\dot\Delta}_{q}\big((\dot\Delta_jg)\circ\psi^{-1}\big)\rangle\big{|}.$$
Therefore  we obtain
$$\big{|}\langle\dot\Delta_j(\dot\Delta_q f\circ\psi),g\rangle\big{|}\leq \|\dot\Delta_{q}f\|_{L^p}
\|\bar{\dot\Delta}_{q}\big((\dot\Delta_jg)\circ\psi^{-1}\big)\|_{L^{\bar p }}.
$$
This implies in view of the first case
\begin{eqnarray*}
\big{|}\langle\dot\Delta_j(\dot\Delta_q f\circ\psi),g\rangle\big{|}&\lesssim& \|\dot\Delta_{q}f\|_{L^p}2^{j-q}\|\nabla\psi^{-1}\|_{L^\infty}\|\dot\Delta_{j}g\|_{L^{\bar p}}\\
&\lesssim&
\|\dot\Delta_{q}f\|_{L^p}2^{j-q}\|\nabla\psi^{-1}\|_{L^\infty}\|g\|_{L^{\bar p}}.
\end{eqnarray*}
Thus we get in view of (\ref{duality1}) the wanted result.
\end{proof}
Next we give a maximum principle estimate for the equation $(TD_{\alpha})$  extending a recent result  due to \cite{C-C} for the partial  case $f=0$. The proof uses  the same idea and  will be briefly described.
\begin{prop}\label{maximum}
Let $v$ be a smooth divergence free vector field and $f$ be a smooth function. We assume that $\theta$ is a smooth solution of the equation
$$
\partial_{t}\theta+v\cdot\nabla \theta+\kappa |\textnormal{D}|^\alpha \theta=f, \quad\textnormal{with}\quad\kappa\geq 0\quad\textnormal{and}\quad
\alpha\in[0,2].
$$
Then for $p\in[1,+\infty]$ we have
$$
\|\theta(t)\|_{L^p}\leq\|\theta(0)\|_{L^p}+\int_{0}^t\|f(\tau)\|_{L^p}d\tau.
$$
\end{prop}
  \begin{proof}
 Let $p\geq 2$, then multiplying the equation by $|\theta|^{p-2}\theta$ and integrating by parts lead to
 $$
\frac{1}{p} \frac{d}{dt}\|\theta(t)\|_{L^p}^p+\kappa\int|\theta|^{p-2}\theta\,|\textnormal{D}|^\alpha\theta dx=\int f |\theta|^{p-2}\theta dx. 
 $$
 On the other hand it is shown in \cite{C-C} that
 $$
 \int|\theta|^{p-2}\theta\,|\textnormal{D}|^\alpha\theta dx\geq 0. 
 $$
 Now using H\"{o}lder's inequality for the right-hand side
 $$
 \int f |\theta|^{p-2}\theta dx\leq \|f\|_{L^p}\|\theta\|_{L^p}^{p-1}. $$
 Thus we obtain
 $$
 \frac{d}{dt}\|\theta(t)\|_{L^p}\leq \|f(t)\|_{L^p}. $$
 We can deduce the result by integrating in time. The case $p\in[1,2[$ can be obtained through  the duality method. 
  \end{proof}
    \begin{Rema}
When this paper was finished we had been informed that similar results were obtained  by Chen et {\it al} \cite{Chen}. In fact they obtained global well-posedness result for small initial data in $\dot B_{p,q}^{s_{c}^p},$ with $p\in[2,\infty[$ and $q\in[1,\infty[.$ For the particular case $q=\infty$ our result is more precise. Indeed, first,  we can extend their result to $p\in[1,\infty]$ and second our smallness condition is given in the space $\dot B_{\infty,1}^{1-\alpha}$ which contains  Besov  spaces $\{\dot{B}_{p,1}^{s_{c}^p}\}_{p\in[1,\infty]}.$
\end{Rema}

\end{document}